\documentclass[dvips,ejs,noinfoline]{imsart}

\RequirePackage[OT1]{fontenc}
\usepackage{amsmath,amsthm,amssymb}
\usepackage[isolatin]{inputenc}
\usepackage{amsfonts,dsfont}
\usepackage{pstricks}
\usepackage{natbib}
\usepackage{color}

\usepackage{graphicx}

\usepackage[colorlinks]{hyperref} 

\doi{10.1214/08-EJS180}
\pubyear{2008}
\volume{2}
\issue{0}
\firstpage{963}
\lastpage{992}

\startlocaldefs
\usepackage{amsthm}
\newtheorem{theorem}{Theorem}[section]
\newtheorem{proposition}[theorem]{Proposition}
\newtheorem{corollary}[theorem]{Corollary}
\newtheorem{lemma}[theorem]{Lemma}

\theoremstyle{remark}
\newtheorem{remark}[theorem]{Remark}
\theoremstyle{definition}

\newtheorem{definition}[theorem]{Definition}

\newcommand{\R}{\mathbb{R}}

\newcommand{\Proba}{\mathbb{P}}
\renewcommand{\P}{\mathbb{P}}

\newcommand{\telque}{\;|\;}
\newcommand{\FDR}{\mathrm{FDR}}

\newcommand{\mbe}{\mathbb{E}}
\newcommand{\mbr}{\mathbb{R}}
\newcommand{\mbp}{\mathbb{P}}

\newcommand{\bp}{\mathbf{p}}
\newcommand{\mtc}{\mathcal}
\newcommand{\mbf}{\mathbf}

\newcommand{\wt}[1]{{\widetilde{#1}}}
\newcommand{\wh}[1]{{\widehat{#1}}}

\newcommand{\ind}[1]{{\mbf{1}\{#1\}}}
\newcommand{\e}[1]{\mbe\brac{#1}}
\newcommand{\ee}[2]{\mbe_{#1}\brac{#2}}
\newcommand{\prob}[1]{\mbp\brac{#1}}
\newcommand{\probb}[2]{\mbp_{#1}\brac{#2}}
\newcommand{\paren}[1]{\left(#1\right)}
\newcommand{\brac}[1]{\left[#1\right]}

\newcommand{\set}[1]{\left\{#1\right\}}
\newcommand{\abs}[1]{\left| #1 \right|}

\newcommand{\al}{\alpha}

\newcommand{\eps}{\varepsilon}

\newcommand{\cH}{{\mtc{H}}}
\newcommand{\cX}{{\mtc{X}}}

\newcommand{\cT}{{\mtc{T}}}
\newcommand{\cN}{{\mtc{N}}}

\endlocaldefs

\begin{document}

\begin{frontmatter}


\title{Two simple sufficient conditions for FDR control}
\runtitle{Sufficient conditions for FDR control}

\begin{aug}
\author{\fnms{Gilles} \snm{Blanchard}\thanksref{t2,t1}\ead[label=e1]{blanchar@first.fraunhofer.de}}
\address{
Fraunhofer-Institut FIRST\\ Kekul\'estrasse 7, 12489 Berlin, Germany \\
\printead{e1}\\
\phantom{E-mail: blanchar@first.fraunhofer.de\ }}
\end{aug}
\and
\begin{aug}
\author{\fnms{Etienne} \snm{Roquain}\corref{}\thanksref{t3,t1}\ead[label=e2]{etienne.roquain@upmc.fr}}
\address{
University of Paris 6, LPMA,\\ 4, Place Jussieu, 75252 Paris cedex 05, France\\
\printead{e2}\\
\phantom{E-mail: etienne.roquain@upmc.fr\ }}

\thankstext{t1}{This work was supported in part by the IST and ICT Programmes of the European
 Community, successively under the PASCAL (IST-2002-506778) and PASCAL2 (ICT-216886) networks of excellence.}
\thankstext{t2}{Part of this work was done while the first author held an invited position at the statistics department of the University of Chicago, which is gratefully acknowledged.}
\thankstext{t3}{Research carried out at the French institute INRA-Jouy and at the Free University of Amsterdam.}

\runauthor{G. Blanchard and E. Roquain}
\end{aug}

\begin{abstract}
 We show that the control of the false discovery rate (FDR) for a
 multiple testing procedure is implied by two coupled simple sufficient
 conditions. The first one, which we call ``self-consistency condition'',
 concerns the algorithm itself, and the second, called ``dependency
 control condition'' is related to the dependency assumptions on the
 $p$-value family. Many standard multiple testing procedures
 are self-consistent (e.g. step-up, step-down or step-up-down procedures), and we prove
 that the dependency control condition can be fulfilled
 when choosing correspondingly appropriate rejection functions, in three classical types
 of dependency: independence, positive dependency (PRDS) and
 unspecified dependency. As a consequence, we recover earlier results
 through simple and unifying proofs while extending their scope to
 several regards: weighted FDR, $p$-value reweighting, new family of
 step-up procedures under unspecified $p$-value dependency and
 adaptive step-up procedures. We give additional examples of other possible applications.
 This framework also allows for defining
 and studying FDR control for multiple testing procedures over a
 continuous, uncountable space of hypotheses.
\end{abstract}



\received{\smonth{2} \syear{2008}}

\end{frontmatter}

\section{Introduction}

A multiple testing procedure is defined as an algorithm taking in
input some (randomly generated) data $X \in \cX$ and returning a set
$R(X)$ of rejected hypotheses, which is a subset of the set $\cH$ of initial candidate
null hypotheses. The false discovery rate (FDR) of the procedure is then defined as the
expected proportion of null hypotheses in $R(X)$ which are in fact true and
thus incorrectly rejected. Following its introduction by \citet{BH1995},
the FDR criterion has emerged recently as a widely used standard
for a majority of applications involving simultaneous testing of a large number of
hypotheses.
It is generally required that a multiple testing procedure $R$
has its FDR bounded by a certain fixed in advance level $\alpha$\,.

Our main point in this work is to show that FDR control is implied by two simple 
conditions. The first one, which we call {\em  self-consistency condition}, requires
that any rejected hypothesis $h\in R(X)$ should have its $p$-value $p_h(X)$ smaller than a
threshold $\Delta_\beta(|R(X)|)$ which itself depends on the volume of rejected hypothesis $|R(X)|$\,,
and on a fixed functional parameter $\beta$\,.
The second one, called {\em dependency control condition}, requires that
for each true null hypothesis $h$\,, the couple of real variables $(U,V)= (p_h,|R(X)|)$
satisfies the inequality (for any $c>0$\,, and the same function $\beta$ as in the first condition):
\begin{equation}
\e{\frac{\ind{U\leq c\beta(V)}}{V}} \leq c\,. \label{equ_probaselfbounding}
\end{equation}
The first condition only concerns how the data is processed to produce the decision, and is hence
purely algorithmic. It can easily be checked
for several classical multiple testing procedures, such as step-down, step-up or more generally step-up-down
procedures. In this condition, the function $\beta$ controls how the threshold
increases with respect to the volume of rejected hypotheses.
In particular, for step-wise procedures, $\beta$ corresponds (up to proportionality constant) to the rejection
function used to cut the curve of ordered $p$-values.
The second condition, on the other hand, is essentially probabilistic in nature.
More precisely, we can show that \eqref{equ_probaselfbounding} can be satisfied
under relatively broad assumptions on the dependency of $(U,V)$\,. In turn,
as will be shown in more detail in the paper,
this implies that the second condition is largely independent of the exact
procedure $R$\,, but rather is related to the dependency assumptions
between the $p$-values.

The two conditions are not independent of each other:
they are coupled through the same functional parameter $\beta$,
appearing in \eqref{equ_probaselfbounding} as well as in the definition of the threshold
$\Delta_\beta$\,. The function $\beta$, called {\em shape function}, is assumed
to be nondecreasing but otherwise arbitrary; if there exists $\beta$ such that the
two corresponding conditions are satisfied, this entails FDR control.

The main advantage of this approach when controlling the FDR is that
it allows us to abstract the particulars of a specific multiple
testing procedure, in order to concentrate on proving the bound
\eqref{equ_probaselfbounding}. This results in  short proofs which in particular do
not resort explicitly to $p$-values reordering.

We then present different types of applications of the result. This approach is first used to show that several
well-known results on FDR control (mainly concerning step-up or step-down procedures based on a linear rejection function)
are recovered in a synthetic way
\citep[e.g., results of][]{BH1995,BH1997,BY2001,Sar2002, GRW2006}.
We also derive the following new results:
\begin{itemize}
 \item[\textbullet] some classical results on step-up procedures are extended to weighted procedures (weighted-FDR and/or $p$-value weighting),
under independence or dependence of the $p$-values;
 \item[\textbullet] a new family of step-up procedures which control the FDR is presented, under unspecified dependencies between the $p$-values;
 \item[\textbullet] we present a simple, exemplary  application of this approach to the problem of {\em adaptive} procedures,
where an estimate of the proportion $\pi_0$ of true null hypotheses in $\cH$ is included in the procedure with the aim of increasing power;
  \item[\textbullet] the case of a continuous space of hypotheses is briefly
 investigated (which can be relevant for instance when the underlying obervation is
 modelled as a stochastic process);
\item[\textbullet] the results of \citet{BL2} and \citet{RS2006} on a specific type of step-down procedure are
 extended to the cases of positive dependencies (under a PRDS-type condition) and unspecified
 dependencies.
\end{itemize}

To put some perspective, let us emphasize here again that the conditions proposed here are only \textit{sufficient}
and certainly not necessary:
naturally, there are many examples of multiple testing procedures that are known to have controlled FDR but do not
satisfy the coupled conditions presented here (including some particular step-up and step-down procedures). The message that
we nevertheless want to convey is that these conditions are able to cover at once an interesting range of
classical existing results on FDR control as well as provide a useful technical tool. It was pointed
out to us that a result similar in spirit to ours will appear in the forthcoming
paper by \citet{FDR2008}; this is discussed in more detail in Section \ref{sec:discuss}.

This paper is organized as follows: in
Section~\ref{sec_prelimin_selfcons}, we introduce the framework, the
two conditions and we prove that taken together, they  imply FDR
control. The self-consistency and
dependency control conditions are then studied separately in Section~\ref{sec:studycond}, leading to specific
assumptions, repectively, on the procedure itself (e.g. step-down, step-up) and on the
dependency between the $p$-values (independence, PRDS, unspecified
dependencies).
The applications summarized above are detailed in Section~\ref{sec:appli}.
Some technical proofs are postponed in the appendix.

\section{Two sufficient conditions for FDR control}\label{sec_prelimin_selfcons}

\subsection{Preliminaries and notations}

Let $(\mtc{X},\mathfrak{X},P)$ be a probability space, with $P$ belonging
to a set or ``model'' $\mathfrak{P}$ of distributions, which can be
parametric or non-parametric.
Formally, a \textit{null hypothesis} is a subset $h \subset \mathfrak{P}$ of distributions on
$(\mathcal{X},\mathfrak{X})$\,. We say that $P$ satisfies
$h$ when $P \in h$\,.

In the multiple testing framework, one is interested in determining
simultaneously whether or not $P$ satisfies
distinct null
hypotheses belonging to a certain set $\cH$\, of candidate hypotheses.  Below, we will always assume
that $\cH$ is at most countable (except specifically in Section
\ref{sec:continuoushypo}, where we mention extensions to continuous
sets of hypotheses).  We denote by $\cH_{0}(P)=\{h\in\cH \telque P \mbox{
satisfies } h\}$ the set of null hypotheses satisfied by $P$, called
the \textit{set of true null hypotheses}. We denote by $\cH_{1}(P)=\cH
\setminus \cH_{0}(P)$ the \textit{set of false null hypotheses} for
$P$\,.

A multiple testing procedure returns a subset $R(x) \subset \cH$ of
rejected hypotheses based on a realization $x$ of a random variable
$X\sim P$\,.
\begin{definition}[Multiple testing procedure]
A \textit{multiple testing procedure} $R$ on $\cH$ is a function $R: x\in \cX  \mapsto R(x) \subset \cH$\,,
such that for any $h \in \cH$\,, the indicator function $\ind{h \in
  R(x)}$ is measurable\,. The hypotheses $h\in R$ are the \textit{rejected null hypotheses} of the procedure $R$.
\end{definition}

We will only consider, as is usually the case, multiple testing
procedures $R$ which can be written as function $R(\bp)$ of a family
of $p$-values $\bp=(p_h, h \in \cH)$\,.  For this, we must assume that for
each null hypothesis $h \in \cH$, there exists a \textit{$p$-value}
function $p_h$\,, defined as a measurable function $p_h: \cX
\rightarrow [0,1]$, such that if $h$ is true, the distribution of
$p_h(X)$ is stochastically lower bounded by a uniform random variable
on $[0,1]$:
\[
\forall P \in \mathfrak{P},\,\qquad \forall h\in \cH_0(P)\,,\forall t\in[0,1]\,, \;\; \probb{X \sim P}{p_h(X) \leq t} \leq t\,.
\]

A {\em type I error} occurs when a true null hypothesis $h$
is wrongly rejected i.e. when $h \in R(x)\cap\cH_0(P)$.  There are several different ways to measure
quantitatively the collective type I error of a multiple testing
procedure.  In this paper, we will exclusively focus on the
false discovery rate (FDR) criterion, introduced by \citet{BH1995} and which has since become a widely used
standard.

The FDR is defined as the averaged proportion of type I errors
 in the set of all the rejected hypotheses.
This ``error proportion'' will be defined in terms of a volume ratio,
and to this end we introduce $\Lambda$\,, a finite positive measure on $\cH$\,.
In the remainder of this paper we will assume such a volume measure has been fixed and
denote, for any subset $S \subset \cH$\,, $|S| = \Lambda(S)$\,.
\begin{definition}[False discovery rate]
Let $R$ be a multiple testing procedure on $\cH$\,. The false discovery rate (FDR) is defined as
\begin{equation}
\label{defFDR}
\FDR(R,P):=\ee{X \sim P}{\frac{|R(X)\cap \cH_0(P)|}{|R(X)|} \ind{|R(X)|>0}}\,.
\end{equation}
\end{definition}
Throughout this paper we will use the following notational convention: whenever there
is an indicator function inside an expectation, this has logical priority over any
other factor appearing in the expectation. What we mean is that if
other factors include expressions that may not be defined (such as
the ratio $\frac{0}{0}$) outside
of the set defined by the indicator, this is safely ignored. In other terms,
any indicator function present implicitly entails that we perform
integration over the corresponding set only. This results in
more compact notation, such as in the above definition.

For the sake of simplifying the exposition, we will (as is usually the
accepted convention) most often drop in the notation a certain number of dependencies,
such as writing $R$ or $p_h$ instead of $R(X)$, $p_h(X)$ and $\cH_0,$ $\cH_1$, $\FDR(R)$
instead of $\cH_0(P)$, $\cH_1(P)$, $\FDR(R,P)$\,. We will
also omit the fact that the probabilities or expectations are
performed with respect to $X \sim P$\,. Generally speaking, we will
implicitly assume that $P$ is fixed, but that all relevant assumptions
and results should in fact hold for any $P \in \mathfrak{P}$\,. For
example, our main goal will be to derive upper bounds on $\FDR(R,P)$
valid for all $P \in \mathfrak{P}$\,; this will be formulated simply
as a bound on $\FDR(R)$\,.

\begin{remark}\label{remark_weightedFDR}(Role of $\Lambda$ and weighted FDR in the finite case)
When the space of hypotheses is finite, the ``standard'' FDR in multiple
testing literature is the one defined using $|.|$ equal to the
counting measure (cardinality) on a finite space and will be referred
to as ``standard $\Lambda$ weigthing''. The notation $|.|$ was kept
here to allow notation compatibility with this case and to alleviate some notational burden.
We stress however that in the case $\cH$ is countably infinite,
the volume measure $\Lambda$ cannot be the cardinality, since we assume it to be finite.

The possibility of using different weights $\Lambda(\set{h})$ for
particular hypotheses $h$ leads to the so-called ``weighted
FDR''. In general, the measure $\Lambda$ represents the relative
importance, or criticality, of committing an error about different hypotheses, and
can be dictated by external constraints.
As discussed by \citet{BH1997} and \citet{BH2006}, controlling
the ``weighted FDR'' can be of interest in some specific
applications. For instance, in the situation where each hypothesis
concerns a whole cluster of voxels in a brain map, it can be relevant to
increase the importance of large discovered clusters when counting the
discoveries in the FDR. Note finally that $\Lambda$ can be rescaled arbitrarily
since only volume ratios matter in the FDR.
\end{remark}

\subsection{Self-consistency, dependency control and the false discovery rate}\label{sec_selfboundingcond}

It is commonly the case that multiple testing procedures are defined as
level sets of the $p$-values:
\begin{equation}
R=\{h\in\cH\telque p_h\leq t\},\label{equ_R}
\end{equation}
where $t$ is a (possibly data-dependent) threshold.
We will be more particularly interested in thresholds that specifically depend on a real parameter $r$ and
possibly on the hypothesis $h$ itself, as introduced in the next definition.
\begin{definition}[Threshold collection]
A \textit{threshold collection} $\Delta$
is a function $${\Delta}: (h,r) \in \cH \times \mbr^+ \mapsto \Delta(h,r) \in \mbr^+,$$
which is nondecreasing in its second variable. A \textit{factorized threshold collection} is a
threshold collection ${\Delta}$ with the particular form: $\forall
(h,r)\in \cH \times \mbr^+,$
\[
\Delta(h,r)=\alpha\pi(h)\beta(r)\,,
\] where
$\pi:\cH\rightarrow [0,1]$ is called the
\textit{weight function} and $\beta:\R^+\rightarrow \R^+$ is a
nondecreasing function called the \textit{shape function}.
Given a threshold collection $\Delta$, the $\Delta$-thresholding-based multiple testing procedure at rejection volume $r$ is defined as
\begin{equation}
L_{\Delta}(r) :=\{h \in \cH \telque p_h \leq \Delta(h,r)\}.\label{equ_thresholdbaseMTP}
\end{equation}
\end{definition}

Let us discuss the role of the parameter $r$\, and proceed to the first of the two announced sufficient conditions.
Remember our goal is to upper bound $\FDR(R)$\,, where the
volume of rejected hypotheses $|R|$ appears as the denominator in the expectation.
Hence, intuitively, whenever this volume gets larger,
we can globally allow more type I errors, and thus take a larger threshold for the $p$-values.
Therefore, the rejection volume parameter $r$ in the definition above should be picked
as an (increasing) function of $|R|$\,.
Formally, this leads to the following ``self-referring'' property:
\begin{definition}[Self-consistency condition]
Given a factorized threshold collection of the form $\Delta(h,r) =
\alpha \pi(h) \beta(r)$\,, a multiple testing procedure $R$ satisfies
the self-consistency condition with respect to the threshold
collection $\Delta$ if the following inclusion holds a.s.:
\begin{equation}
\label{SCC}
R\subset L_\Delta(|R|).\tag{{\bf SC}$(\alpha,\pi,\beta)$}
\end{equation}
\end{definition}

Next, we introduce the following probabilistic condition on two dependent real variables:
\begin{definition}[Dependency control condition]
Let $\beta: \mbr^+ \rightarrow \mbr^+$ be a nondecreasing function.  A
couple $(U,V)$\, of (possibly dependent) nonnegative real random variables is said to
satisfy the dependency control condition with shape function $\beta$
if the following inequalities hold:
\begin{equation}
\label{DCC}
\forall c>0,\:\:\:\:\e{\frac{\ind{U \leq c \beta(V)}}{V}} \leq c\,. \tag{{\bf DC}$(\beta)$}
\end{equation}
\end{definition}

The following elementary but fundamental result is the main cornerstone
linking the FDR control to conditions {\bf SC} and {\bf DC}.
\begin{proposition}
\label{mlem}
Let $\beta: \mbr^+ \rightarrow \mbr^+$ be a (nondecreasing) shape function,
$\pi: \cH \rightarrow [0,1]$ a weight function and $\alpha$ a positive number\,.
Assume that the multiple testing procedure $R$ is such that:

\smallskip

\noindent (i) the self-consistency condition \ref{SCC} is satisfied\,;

\noindent (ii) for any $h \in \cH_0$\, the couple $(p_h,|R|)$ satisfies \ref{DCC}.

\smallskip

\noindent Then $\FDR(R) \leq \al \Pi(\cH_0)$\,, where $d\Pi = \pi d\Lambda$\,, i.e., $\Pi(\cH_0) := \sum_{h \in \cH_0} \Lambda(\set{h}) \pi(h)$\,.
\end{proposition}
\begin{proof}
From \eqref{defFDR},
\begin{align*}
\FDR(R)= \e{\frac{|R\cap \cH_0|}{|R|} \ind{|R|>0}} & = \sum_{h\in\cH_0} \Lambda(\set{h}) \e{\frac{\ind{h\in R}}{|R|}}
\\
&\leq\sum_{h\in\cH_0} \Lambda(\set{h})
  \e{\frac{\ind{p_h\leq \alpha\pi(h)\beta(|R|)}}{|R|}}\\
  &\leq  \alpha \sum_{h\in\cH_0} \Lambda(\set{h}) \pi(h) ,
\end{align*}
where we have used successively conditions (i) and (ii) for the two
above inequalities.
\end{proof}

Let us point out the important difference in nature between the two sufficient conditions:
for a fixed shape function $\beta$, the self-consistency condition (i) concerns
only the algorithm itself (and not the random structure of the problem).
On the other hand, the dependency control condition (ii) seems to
involve both the algorithm and the statistical nature of
the problem. However, we will show below in Section \ref{sec:depcond}
that this latter condition can be checked
under a weak, general and quite natural assumption on the algorithm itself (namely that
$|R(\bp)|$ is nonincreasing function of the $p$-values), and primarily depends on the
dependency structure of the $p$-values. (Moreover, in the case of arbitrary dependencies,
we will consider a special family of $\beta$s which satisfy the condition without any assumptions on the algorithm.)
Hence, the interest of the above proposition is that it effectively
separates the problem of FDR control into a purely {\em algorithmic} and an (almost) purely
{\em probabilistic} sufficient condition. The link between the two conditions is the common
shape function $\beta$\,: the dependency assumptions between the $p$-values will
determine for which shape function the condition \ref{DCC} is valid; in turn, this
will impose constraints on the algorithm through condition \ref{SCC}.

\begin{remark}\label{rem_pi} 
(Role of $\pi$ and $p$-value weighting in the finite case)
To understand intuitively the role of the weight function $\pi$\,, assume $\cH$ is of finite cardinality $m$ and
take for simplification $\beta(r)=1$ for now. Consider the corresponding testing procedure $L_\Delta$:
the rejected hypotheses are those for which $p'_h := p_h / (m \pi(h)) \leq \alpha/m$\,,
where $p'_h$ is the {\em weighted $p$-value} of $h$\,. If $\pi(h)$ is constant equal to $1/m$\,, we have $p'_h=p_h$ and
the above is just Bonferroni's procedure, which has family-wise error rate (FWER) controlled by $\alpha$\,.
If $\pi(h)$ is, more generally, an arbitrary probability distribution on $\cH$\,,
the above is a weighted Bonferroni's procedure and has also FWER less than $\alpha$ \citep[see, e.g.,][]{WR2006}. In this example, $\pi$ represents
the relative importance, or weight of evidence, that is given {\em a priori} to $p$-values, and thus plays the role
of a prior that can be fixed arbitrarily by the user. Its role in the control of FDR is very similar; the use
of weighted $p$-values for FDR control has been proposed earlier, for example by \citet{GRW2006}.
When $\cH$ is of finite cardinality $m$, we will refer to the choice $\pi(h)\equiv 1/m$\,, in conjunction with $\Lambda$ being the cardinality measure, as
the ``standard $\Lambda-\pi$ weighting''.

More generally, following Proposition~\ref{mlem}, control of the FDR at level $\alpha$ is ensured as soon as
the weight function $\pi$ is chosen as a probability density with respect to $\Lambda$
(i.e. $\sum_{h\in\cH} \Lambda(\set{h})\pi(h)=1$).
When $\cH$ is of finite cardinality $m$ and with the ``standard $\Lambda-\pi$ weighting'' defined above,
we obtain
$\FDR\leq \alpha m_0/m \leq \alpha$
(where $m_0$ denotes the number of true null hypotheses).
\end{remark}

\begin{remark}\label{rem_mlem} 
Proposition~\ref{mlem} can be readily extended to the case
where we use different volume measures for the numerator and
denominator of the $\FDR$\,. However, since it is not clear to us whether
such an extension would be of practical interest, we choose in this paper to
deal only with a single volume measure.
\end{remark}

\section{Study of the two sufficient conditions}\label{sec:studycond}

In this section, we give a closer look to conditions \ref{SCC} and \ref{DCC},
and study typical situations where they are statisfied.

\subsection{Self-consistency condition and step-up procedures}
\label{section_stepup}

The main examples of self-consistent procedures are step-up procedures. In fact,
for a fixed choice of parameters $(\alpha,\beta,\pi)$\,, step-up procedures output
the largest set of rejected hypotheses such that \ref{SCC} is satisfied, and are in this
sense optimal with respect to that condition. Here, we define step-up procedures
by this characterizing property, thus avoiding the usual definition using the reordering
of the $p$-values.

\begin{definition}[Step-up procedure]\label{def_stepup_chapself}
Let $\Delta$ be a factorized threshold collection of the form $\Delta(h,r) =
\alpha \pi(h) \beta(r)$\,.
The \textit{step-up multiple testing procedure} $R$ associated to
$\Delta$\,,
is given
by either of the following equivalent definitions:
\begin{align*}
(i) \:\:\:& R = L_\Delta(\hat{r})\,,\mbox{ where }
\hat{r}:=\max\{r \geq 0 \telque |L_\Delta(r)|\geq r\}\,\\
(ii) \:\:\:&
R=\bigcup\big\{A \subset \cH \telque A \mbox{ satisfies
\ref{SCC} }\big\}\,.
\end{align*}
Additionally, $\wh{r}$ satisfies $|L_\Delta(\wh{r})| =\wh{r}$\,; equivalently,
the step-up procedure $R$ satisfies \ref{SCC} with equality.
\end{definition}

\begin{proof}[Proof of the equivalence between $(i)$ and $(ii)$]
Note that, since $\Delta$ is assumed
  to be nondecreasing in its second variable, $L_\Delta(r)$ is a nondecreasing set
  as a function of $r\geq 0$\,. Therefore, $\abs{L_\Delta(r)}$ is a nondecreasing function
of $r$ and the supremum appearing
in $(i)$ is a maximum i.e. $|L_\Delta(\hat{r})|\geq \hat{r}$\,.
It is easy to see that $|L_\Delta(\hat{r})| = \hat{r}$ because this would
otherwise contradict the definition of $\wh{r}$\,.
Hence $L_\Delta(\hat{r}) =
L_\Delta(|L_\Delta(\hat{r})|)$\,,
so $L_\Delta(\hat{r})$ satisfies \ref{SCC} (with equality) and is included in the set union appearing in
$(ii)$. Conversely, for any set $A$ satisfying $A \subset
L_\Delta(|A|)$\,, we have $|L_\Delta(|A|)| \geq |A|$\,, so that $|A|
\leq \hat{r}$ and $A \subset L_\Delta(\hat{r})$\,.
\end{proof}

\begin{figure}[t]
\resizebox{.9\textwidth}{!}{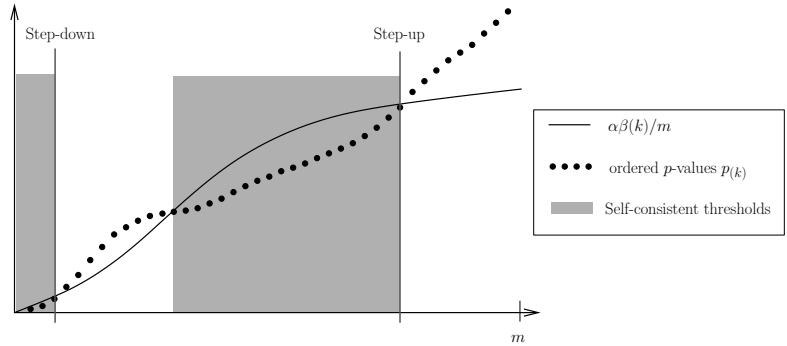}
\caption{\label{fignoob} Pictorial representation of the step-up (and step-down) thresholds, and (in grey)
of all thresholds $r\in\{1,\ldots,m\}$ for which $L_\Delta(r)$ satisfies the self-consistency condition.
The $p$-values and the rejection function represented here  have been picked arbitrarily and in a deliberately exaggerated fashion
in order to illustrate the different procedures; they are not meant to represent a realistic data or model.
This picture corresponds to the standard $\Lambda$-$\pi$ weighting only.
}
\end{figure}

When $\cH$ is finite of cardinal $m$ endowed with the standard $\Lambda$-weighting $\Lambda(\cdot)=\mathrm{Card}(\cdot)$\,,
Definition~\ref{def_stepup_chapself} is equivalent to the
classical definition of a step-up procedure, based on reordering the $p$-values:
for any $h\in\cH$, denote by $p'_h:=p_h/(m\pi(h))$ the
\textit{weighted $p$-value} of $h$
(in the case $\pi(h)=0$, we put $p'_h=+\infty$ if $p_h>0$ and $p'_h=0$ if $p_h=0$), and consider the
ordered weighted $p$-values
$$p'_{(1)}\leq p'_{(2)}\leq \dots \leq p'_{(m)}.$$
Since $L_{\Delta}(r)=\{h\in\cH\telque p'_h\leq \alpha\beta(r)/m\}$,
the condition $|L_{\Delta}(r)|\geq r$ is equivalent to $p'_{(r)}\leq
\alpha \beta(r)/m$\,. Hence, the step-up procedure associated to
$\Delta$ defined in Definition \ref{def_stepup_chapself} rejects all
the $\hat{r}$ smallest weighted $p$-values, where $\hat{r}$
corresponds to the ``last right crossing'' point between the ordered
weighted $p$-values $p'_{(\cdot)}$ and the scaled shape function $\alpha
\beta(\cdot)/m$:
\begin{eqnarray}
\hat{r}&=&\max\big\{r\in\{0,\dots,m\}\telque p'_{(r)}\leq \alpha \beta(r)/m\big\}, \nonumber
\end{eqnarray}
with $p'_{(0)}:=0$; see Figure \ref{fignoob} for an illustration. For the standard $\pi$-weighting $\pi(h)=1/m$, the
weighted $p$-values are simply the $p$-values. In particular:
\begin{itemize}
\item The step-up procedure associated to the linear shape function $\beta(r)=r$ is the well-known linear
step-up procedure of \citet{BH1995}.
\item The step-up procedure associated to the linear shape function
$\beta(r)=r \paren{\sum_{i=1}^m \frac{1}{i}}^{-1}$ is the
distribution-free linear step-up procedure of Benjamini and Yekutieli (\citeyear{BY2001}).
\end{itemize}

Finally, let us point out that {\em step-down} and more generally {\em step-up-down}
procedures are also self-consistent. The latter class of step-wise
procedures have been introduced by \citet{TLD1998}, and contains
step-up and step-down procedures as particular cases.
These procedures select in a certain way
among the ``crossing points''  between the $p$-value function and some fixed rejection function
(for example, on Figure~\ref{fignoob}, there are only two non-zero crossing points to choose from).
More formally, and under arbitrary weighting, given a parameter $\lambda\in [0,|\cH|]$, the step-up-down procedure with threshold
collection $\Delta$ and of order $\lambda$ is defined as
$L_\Delta(\widehat{r}_\lambda)$, where either $\wh{r}_\lambda:=\max\{r
\geq \lambda \telque \forall r', \lambda \leq r'\leq r,
|L_\Delta(r')|\geq r'\}$ if $|L_\Delta(\lambda)|\geq \lambda$\,; or
$\wh{r}_\lambda:=\max\{r < \lambda \telque |L_\Delta(r)|\geq r\}$
otherwise.  In words, assuming the standard weighting case and $\lambda$ an integer,
if $p_{(\lambda)}$ is smaller than the rejection function at $\lambda$\,,
the closest crossing point to the right of $\lambda$ is picked, otherwise the closest crossing point to the left.
In particular, the step-up-down procedure of order
$\lambda=|\cH|$ is simply the step-up procedure (based on the same
threshold collection). The case $\lambda=0$ is the \textit{step-down}
procedure.
Although generalized step-up-down procedures are not maximal with respect to condition ${\bf SC}$
like the plain step-up,
the fact that they still satisfy that condition is worth noticing.

\subsection{Dependency control condition}
\label{sec:depcond}

In this section, we show that condition (ii) of Proposition
\ref{mlem} holds under different types of assumptions on the
dependency of the $p$-values. We will follow
the different types of dependencies considered by \citet{BY2001},
namely independent, positive dependency under the PRDS
condition and arbitrarily dependent $p$-values.
In each case, we have to prove \ref{DCC}
for specific conditions on the variables $(U,V)$\,,
resulting in specific choices for the shape function $\beta$\,.

We start the section with a probabilistic lemma collecting the technical tools used
to deal with each situation.
\begin{lemma}
\label{UVlemma}
Let $(U,V)$ be a couple of nonnegative random variables such that $U$ is stochastically
lower bounded by a uniform variable on $[0,1]$\,, i.e. $\forall t\in [0,1], \P(U\leq t)\leq t$\,. Then the dependency control condition
\ref{DCC} is satisfied by $(U,V)$ under any of the following situations:

\smallskip

\noindent(i) $\beta(x)=x$ and $V=g(U)$\,, where $g:\mathbb{R}^+ \rightarrow \mathbb{R}^+$ is a nonincreasing function.

\smallskip

\noindent(ii) $\beta(x)=x$ and the conditional distribution of $V$ given ${U \leq u}$ is stochastically
decreasing in $u$, that is,
\begin{equation}
\label{UV-nondec}
\text{ for any } r\geq 0\,, \text{ the function }u\mapsto \Proba(V <
r\telque U \leq u) \text{ is nondecreasing}\,.
\end{equation}

\smallskip

\noindent(iii) The shape function is of the form
\begin{equation}
\label{condbeta_chapself}
\beta_\nu(r) = \int_0^r x d\nu(x)\,,
\end{equation}
where $\nu$ is an arbitrary probability distribution on $(0,\infty)$\,, and $V$ is arbitrary.
\end{lemma}
The proof is found in appendix.
Note that there is some redundancy in the lemma since (i) is a particular case of (ii), but
this subcase has a particularly simple proof and is of self interest because it corresponds to the case
of independent $p$-values (as will be detailed below).

We now apply this result to prove
that for any $h\in\cH_0$, the couple of variables $(p_h,|R|)$ satisfies \ref{DCC},
under the different dependency assumptions on the $p$-values, and for the correspondingly
appropriate functions $\beta$ given by the lemma.
The only additional assumption we will make on the procedure $R$ itself is that
it has nonincreasing volume as a function of the $p$-values (and this assumption will
not be required in the case of arbitrarily dependent $p$-values).

\subsubsection{Independent case}

\begin{proposition}
\label{propCCindep}
Assume that the
collection of $p$-values $\bp=(p_h,h\in\cH)$
forms an independent family of random variables.
Let $R(\bp)$ be a multiple testing procedure
such that $|R(\bp)|$ is nonincreasing in each $p$-value $p_h$ such that
$h\in\cH_0$\,. For any $h \in \cH$, denote $\bp_{-h}$ the collection of $p$-values
 $(p_{g} : g \in \cH, g\neq h)$\,.

Then for any $h \in \cH_0$ and for the linear shape function $\beta(x)=x$\,, the couple of variables $(p_h,|R|)$ satisfies
\ref{DCC}, in which the expectation is taken conditionally to the $p$-values of $\bp_{-h}$. As a consequence, it also satisfies \ref{DCC} unconditionally.
\end{proposition}

\begin{proof}
By the independence assumption,
the distribution of $U=p_h$ conditionally to $\mbf{p}_{-h}$ is identical
to its marginal and therefore stochastically lower bounded by a
uniform distribution. The value of $\bp_{-h}$ being held fixed, $|R(\bp)| = |R((\bp_{-h},p_h))|$
can be written as a nonincreasing function $g$ of $p_h$ by the assumption
on $R$\,. We conclude by part (i) of Lemma \ref{UVlemma}.
\end{proof}

\begin{remark}
Note that Proposition \ref{propCCindep} is still valid under the
slightly weaker assumption that for all $h\in\cH_0$, $p_{h}$ is
independent of the family $(p_{g},g\neq h)$ (in particular, the $p$-values of
$(p_{h},h\in\cH_1)$ need not be mutually independent).
\end{remark}

\subsubsection{Positive dependencies (PRDS)}
\label{prds_dep}

From point (ii) of Lemma \ref{UVlemma}, each couple $(p_h,|R|)$ satisfies \ref{DCC} with $\beta(x)=x$
under the following condition (weaker than independence):
\begin{equation}
\label{R-nondec}
\mbox{ for any }r\geq 0\,, \mbox{ the function }u\mapsto \Proba(|R| <
r\telque p_h \leq u) \mbox{ is nondecreasing}\,.
\end{equation}

Following \citet{BY2001},
we state a dependency condition ensuring that $(p_h,|R|)$ satisfies (\ref{R-nondec}).
For this, we recall the definition of positive regression dependency
on each one from a subset (PRDS) (introduced by \citealp{BY2001}, where
its relationship to other notions of positive dependency is also
discussed).  Remember that a subset $D\subset[0,1]^\cH$ is called
\textit{nondecreasing} if for all $\mbf{z},\mbf{z}' \in [0,1]^\cH$
such that $\mbf{z} \leq \mbf{z}'$ (i.e. $\forall h\in\cH, z_h\leq
z'_h$), we have $\mbf{z} \in D \Rightarrow \mbf{z}' \in D$\,.

\begin{definition}\label{def_PRDS}
For $\cH'$ a subset of $\cH$\,, the $p$-values of
$\mathbf{p}=(p_h,h\in \cH)$ are said to be \textit{positively
  regressively dependent on each one from} $\cH'$ (denoted in short by
\textit{PRDS on $\cH'$}), if for any $h\in\cH'$\,,
for any measurable nondecreasing set $D\subset[0,1]^{\cH}$\,, the function $ u
\mapsto \Proba(\mbf{p}\in D\telque p_h=u)$\,
is nondecreasing.
\end{definition}

We can now state the following proposition:
\begin{proposition}\label{cor_PRDSimpliesfdrcontrol}
Suppose that the $p$-values of $\mathbf{p}=(p_h,h\in\cH)$ are PRDS on $\cH_0$\,,
and consider a multiple testing procedure $R$
such that $|R(\bp)|$ is nonincreasing in each $p$-value.
Then for any $h \in \cH_0$\,, the couple of variables $(p_h,|R|)$ satisfies
\ref{DCC} for the linear shape function $\beta(x)=x$\,.
\end{proposition}
\begin{proof}
We merely check that condition \eqref{R-nondec} is satisfied.
For any fixed $r\geq 0$\,, put $D = \set{\mbf{z}\in[0,1]^{\cH}\telque |R(\mbf{z})| < r}$\,. It is clear
from the assumptions on $R$ that $D$ is a nondecreasing measurable set. Then
by elementary considerations, the PRDS condition (applied using the set $D$ defined above)
implies \eqref{R-nondec}. The latter argument was also used by \citet{BY2001} with a reference
to \citet{Lehm1966}. We provide here a succinct proof of this fact in the interest of remaining self-contained.

 Under the PRDS condition, for all $u\leq u'$\,,
 putting $\gamma = \prob{p_h \leq u \telque p_h \leq u'}$\,,
 \begin{align*}\
 \prob{ \mbf{p} \in D \;|\; p_h \leq u'}
 & = \e{ \prob{ \mbf{p} \in D \;|\; p_h} \;|\; p_h \leq u'}\\
 & = \gamma \e{ \prob{\mbf{p} \in D \;|\;p_h} \;|\; p_h \leq u} \\
& \;\;\;\; + (1-\gamma) \e{ \prob{\mbf{p} \in D \;|\;p_h} \;|\; u < p_h \leq u'}\\
 & \geq \e{ \prob{\mbf{p} \in D \;|\;p_h} \;|\; p_h \leq u}
 = \prob{ \mbf{p} \in D \;|\; p_h \leq u}\,,
 \end{align*}
where we have used the definition of PRDS for the last inequality.
\end{proof}

\subsubsection{Unspecified dependencies}
\label{sssec:unsp}

We now consider a totally generic
 setting with no assumption on the
dependency structure between the $p$-values nor on the structure of
the multiple testing procedure $R$\,. Using point (iii) of Lemma \ref{UVlemma},
we obtain immediately the following result:
\begin{proposition}\label{lemmaCCdep}
Let $\beta_\nu$ be a shape function of the form \eqref{condbeta_chapself}.
Then for any $h \in \cH_0$\,, the couple of variables $(p_h,|R|)$ satisfies
\ref{DCC}, for any multiple testing procedure $R$\,.
\end{proposition}
Note that a shape function of the form \eqref{condbeta_chapself}
must satisfy  $\beta_\nu(r)\leq r$\,, with strict inequality except for at most one
point beside zero (some examples will be discussed below in Section \ref{ssec:unspec}).
Therefore, the price to pay here
is a more conservative dependency control inequality, in turn
resulting in a more restrictive self-consistency condition
when using this shape function.
This form of shape function was initially
introduced by \citet{BF2007}, where some ties were exposed between multiple testing
and statistical learning theory.

\section{Applications}\label{sec:appli}

\subsection{The linear step-up procedure with $\Lambda-\pi$ weighting }\label{applistep-up}

We have seen earlier in Section \ref{section_stepup} that step-up procedures satisfy the self-consistency condition.
Furthermore, is is easy to see that step-up procedures are nonincreasing as a function of the $p$-values.
Using this in conjunction with Proposition~\ref{propCCindep} (resp. Proposition~\ref{cor_PRDSimpliesfdrcontrol})
and Proposition~\ref{mlem},
we obtain the following result for the ($\Lambda$-weighted) FDR control of
the ($\pi$-weighted) linear
step-up procedure -- that is, the step-up procedure associated to the
threshold collection $\Delta(h,r)=\alpha\pi(h)r$\,.
\begin{theorem} \label{th_stepup}
For any finite and positive volume measure $\Lambda$\,,
the ($\pi$-weighted) linear step-up procedure $R$ has its ($\Lambda$-weighted) FDR
upper bounded by $\Pi(\cH_0)\alpha$\,, where
$\Pi(\cH_0):=\sum_{h\in\cH_0}\Lambda(\set{h}) \pi(h)$, in either of the following
cases:
\begin{itemize}
\item the $p$-values of $\mathbf{p}=(p_h,h\in \cH)$ are independent.
\item the $p$-values of $\mathbf{p}=(p_h,h\in \cH)$ are PRDS on $\cH_0$.
\end{itemize}
\end{theorem}
Again, the statement is redundant since independence is a particular case of PRDS,
and we just wanted to recall that the treatment of the independent case is particularly simple.
This theorem essentially recovers and unifies some known results concerning
particular cases: the two points of the theorem were respectively
proved by \citet{BH1995} and \citet{BY2001}, with a uniform $\pi$, and
$\Lambda$ the cardinality measure.  For a general volume measure $\Lambda$ and a
uniform $\pi$, the above result in the independent case was proved by
\citet{BH1997}.  A proof with a general $\pi$, $\Lambda$ the cardinality measure
and in the independent case was investigated by \citet{GRW2006}.

The interest of the present framework is to  allow for a general and unified
version of these results with a concise proof (avoiding in particular
to consider explicitly $p$-value reordering). We distinguish clearly between
the two different ways to obtain ``weighted'' versions of step-up procedures,
by changing respectively the choice of the volume measure $\Lambda$
or the weight function $\pi$.
Both types of weighting are of interest and of different nature; using weighted $p$-values can
have a large impact on power (\citealp{GRW2006,RW2008}; see also above Remark~\ref{rem_pi}), while using a
volume $\Lambda$ different from the cardinality measure can be of relevance
for some application cases (see Benjamini and Hochberg, \citeyear{BH1997}; \citealp{BH2006}; and Remark \ref{remark_weightedFDR} above).
Up to our knowledge, the two types of weighting had not been considered simultaneouly before;
in particular and
as  noticed earlier (see Remark~\ref{rem_pi}), in order to ensure FDR control at level $\alpha$\, under an
arbitrary volume measure $\Lambda$\,, the appropriate choice for a
weight function $\pi$ is to take a density function with respect to
$\Lambda$\,.\looseness=-1

\subsection{An extended family of step-up procedures under unspecified dependencies}
\label{ssec:unspec}

Similarly, in the case where the $p$-values have unspecified dependencies, we use
Proposition \ref{lemmaCCdep} instead of Proposition \ref{cor_PRDSimpliesfdrcontrol}
to derive the following theorem:
\begin{theorem}\label{cor_marteau}
Consider $R$ the step-up procedure associated to the factorized
threshold collection $\Delta(h,r)=\alpha \pi(h) \beta_\nu(r)$, where the
shape function $\beta_\nu$ can be written in the form \eqref{condbeta_chapself}.
Then $R$ has its ($\Lambda$-weighted) FDR controlled at level $\Pi(\cH_0)\alpha$\,.
\end{theorem}

Theorem \ref{cor_marteau} can be seen as an extension to the FDR
of a celebrated inequality due to
\citet{Hom1983} for the family-wise error rate (FWER), which has been
widely used in the multiple testing literature
\citep[see, e.g.,][]{LR2005,RS2006,RS2006b}.  Namely, when $\nu$ has its
support in $\{1,\dots,m\}$ and $\cH = \cH_0$\,, the above result recovers
Hommel's inequality. Note that the latter special  case corresponds to a
``weak control'', where we assume that all null hypotheses are true; in
this situation the FDR is equal to the FWER. Note also that Theorem \ref{cor_marteau} generalizes without
modification to a possibly continuous hypothesis space, as will be
mentioned in Section \ref{sec:continuoushypo}. The result of Theorem~\ref{cor_marteau} initially appeared
in a paper of \citet{BF2007}, in a somewhat different setting.

\subsubsection{Discussion of the family of new shape functions}

Theorem \ref{cor_marteau} establishes that, under arbitrary
dependencies between the $p$-values, there exists a family of step-up
procedures with
controlled false discovery rate. This family is parametrized
by the free choice
of a distribution $\nu$ on the positive real line, which determines
the shape function $\beta_\nu$\,.

In the remaining of Section~\ref{ssec:unspec}, we assume $\cH$ to be finite of cardinal $m$\,, endowed
with the standard $\Lambda$ weighting, i.e., the counting measure.
In this situation, let us first remark that it is always preferable to choose $\nu$ with support in $\set{1,\ldots,m}$\,. To see this, notice
that only the values of $\beta$ at integer values $k,  1\leq k \leq m$\, matter for the output of
the algorithm. Replacing an arbitrary distribution $\nu$ by the discretized distribution
$\nu'(\set{k}) = \nu((k-1,k])$ for $k<m$ and $\nu'(\set{m}) = \nu((m-1,+\infty))$
results in a shape function $\beta'$ which is larger than $\beta$ on the relevant integer range,
hence the associated step-up procedure is more powerful. This discretization operation
will however generally result in minute improvements only; sometimes continuous distributions
can be easier to handle and avoid cumbersomeness in theoretical considerations.

Here are some simple possible choices for (discrete) $\nu$\, based on power functions $\nu(\{k\})\propto k^\gamma$\,, $\gamma \in \set{-1,0,1}$\,:
\begin{itemize}
\item $\nu(\set{k})=\gamma_m^{-1}k^{-1}$ for $k\in\{1,\dots,m\}$ with the normalization
constant $\gamma_m=\sum_{1\leq i\leq m}\frac{1}{i}$. This yields
$\beta(r)= \gamma_m^{-1}r\,,$ and we recover the distribution-free procedure of
\citet{BY2001}.
\item $\nu$ is the uniform
on $\{1,\dots,m\}$, giving rise to the quadratic shape function
$\beta(r)=r(r+1)/2m\,.$ The obtained step-up procedure was proposed by \citet{Sar2006}.
\item $\nu(\set{k})=2k/(m(m+1))$ for $k\in\{1,\dots,m\}$
leads to $\beta(r)=r(r+1)(2r+1)/(3m(m+1))\,.$
\end{itemize}
\eject

\begin{figure}[t!]
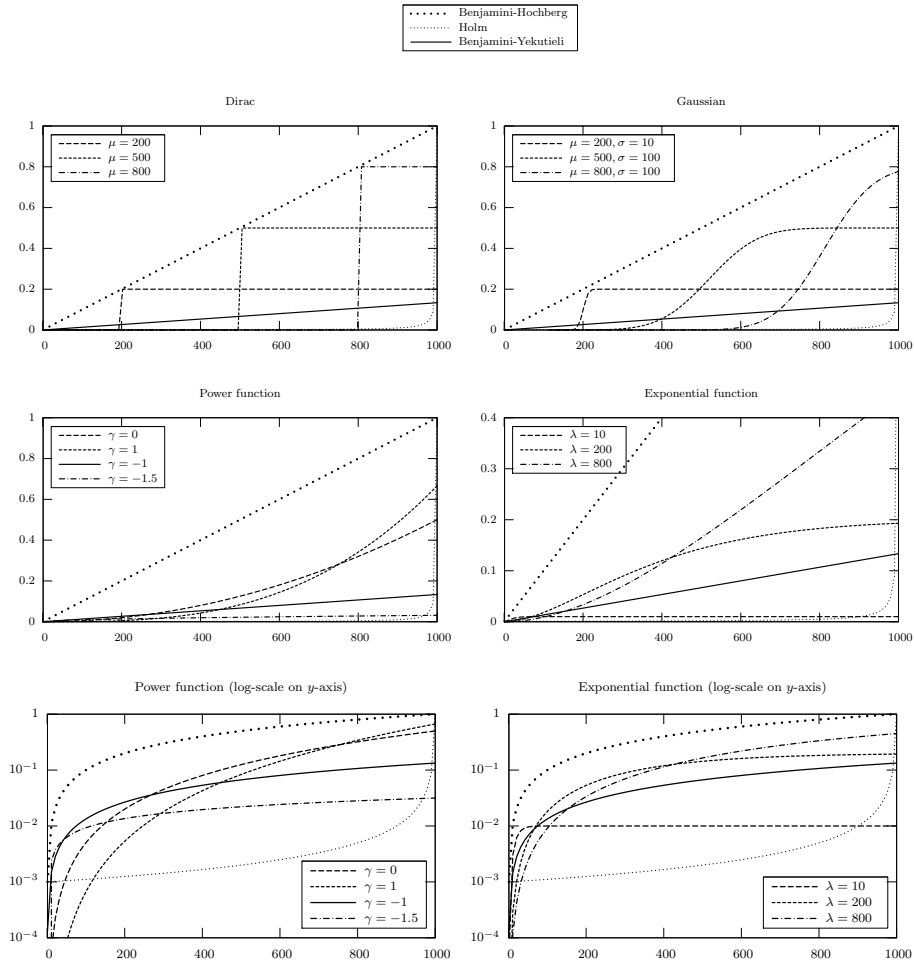

\vspace{30pt}
\begin{center}
\makebox[.3\textwidth][r]{\smash{\resizebox{.51\textwidth}{!}{\input{commonlegend.tex}}}}\\
\begin{tabular}{cc}
\hspace{-7mm} \resizebox{.52\textwidth}{!}{\input{gpseuildirac.tex}} &
\hspace{-7mm} \resizebox{.52\textwidth}{!}{\input{gpseuilgauss.tex}}\\
\hspace{-7mm} \resizebox{.52\textwidth}{!}{\input{gpseuilpower.tex}} &
\hspace{-7mm} \resizebox{.52\textwidth}{!}{\input{gpseuilexp.tex}}\\
\hspace{-13mm} \resizebox{.57\textwidth}{!}{\input{gpseuilpower_log.tex}} &
\hspace{-13mm} \resizebox{.57\textwidth}{!}{\input{gpseuilexp_log.tex}}\\
\end{tabular}
\end{center}\vspace{-2mm}
\caption{\label{fig_choicenu} For the standard $\Lambda$-weighting and $m=1000$ hypotheses, this figure shows
several (normalized) shape functions $m^{-1}\beta$ associated to different
distributions $\nu$ on $\R^+$ (according to expression \eqref{condbeta_chapself}):
\textit{Dirac distribution}: $\nu=\delta_{\mu},$ with $\mu>0$.
\textit{(Truncated-) Gaussian distribution}: $\nu$ is
  the distribution of $\max(X,1)$, where $X \sim \cN(\mu,\sigma^2)$\,.
\textit{Power distribution}: $d\nu(r)=r^\gamma\ind{r\in[1,m]}dr/\int_{1}^{m}u^{\gamma}du$,
  $\gamma\in\R$.
\textit{(Truncated-) Exponential distribution}: $d\nu(r)=(1/\lambda)\exp(-r/\lambda)\ind{r\in[0,m]}dr,$ with
  $\lambda>0$. On each graph, for comparison purposes we added the threshold function for Holm's step-down
$m^{-1}\beta(x) = 1/(m-x+1)$\,, (small dots),
and the linear thresholds $\beta(x)=x$ (large dots) and
$\beta(x) = (\sum_{i\leq m} i^{-1})^{-1} x$\, (solid -- also corresponding to the power distribution with $\gamma=-1$),
corresponding to the standard linear step-up and to the distribution-free
linear step-up of \citet{BY2001}, respectively.}
\vspace{-12pt}
\end{figure}

On Figure \ref{fig_choicenu}, we plotted the shape functions corresponding
to different choices of distributions $\nu$ (which are actually  continuous, \textit{i.e.}, without applying the
discretization procedure mentioned above).
It is clear that the choice of $\nu$ has a large impact on the final number of rejections of the
procedure. However, since no shape function uniformly dominates the others,
there is no universally optimal choice of $\nu$: the
respective performances of these different procedures will depend on the exact distribution $P$\,,
and in particular on the number of non-true hypotheses.

We like to think of $\nu$ as a kind of ``prior'' on the possible volumes of rejected hypotheses.
If we expect to have only a few rejected hypotheses, $\nu$ should be concentrated on small values,
and more spread out if we expect a significant rejection proportion.
This intuition is in accordance with a case of equality in Hommel's inequality established
by \citet[Lemma 3.1 (ii)\,]{LR2005}. In the situation studied there (a specifically crafted
distribution $P$\,), it can be checked that the distribution
of the cardinality of the step-up procedure $R$ using the shape function $\beta_\nu$\,,
conditionally to $R\neq \emptyset$\,, is precisely $\nu$\, in our notation, while $\FDR(R)$ is exactly $\alpha$\,.

As mentioned previously in Section \ref{sssec:unsp}, for any choice of
$\nu$\,, the shape function $\beta_\nu$ is always upper bounded by the
linear shape function $\beta(x)=x$\,. The only cases of equality
are attained if $\nu$   is equal to a Dirac measure $\delta_{x_0}$ in a
point $x_0 \in
\{1,\ldots,m\}$\,: in this case $\beta_{\delta_{x_0}}(x_0) = x_0$ but
$\beta_{\delta_{x_0}}(x)<x$ for any $x\neq x_0$\,. Therefore, these procedures
always reject less (or at most as many) hypotheses than the linear
step-up. Admittedly, this probably limits the practical implications
of this result, as we expect practitioners to prefer using the
standard linear step-up even if the theoretical conditions for its
validity cannot be formally checked in general.  Additional conservativeness is the
``price to pay'' for validity under arbitrary dependencies, although
the above result shows that
one has, so to say, the choice in the way this price is to be paid.

Finally, from the examples of shape functions drawn on Figure
 \ref{fig_choicenu}, the shape functions based on {\em exponential} distributions $\nu$\, seem particularly
interesting; they appear to exhibit a qualitatively diverse range of
 possible shape functions, offering more flexibility than
the Benjamini-Yekutieli procedure while not being as committed as the
 Dirac distributions to a
 specific prior belief on the number of rejected hypotheses.

\subsubsection{Comparison to Bonferroni's and Holm's procedures}

Observe that Bonferroni's procedure also belongs to the family presented here
(taking $\nu=\delta_1$) -- in the sense that a single-step procedure using a fixed threshold can
be technically considered as a step-up procedure.
It is well-known, however, that its control on type I
error is much stronger than bounded FDR, namely bounded FWER. To this extent, it is worth
considering the question of whether other rejections functions in the family -- for which only the FDR is controlled -- are of interest at all
As remarked earlier,
no shape function in the family can \textit{uniformly} dominate the others, and consequently
there exist particular situations where Bonferroni's procedure can be more powerful (i.e. reject more hypotheses)
than other members of the family.
However, this case appears
only
when there is indeed a very small number of rejections
(i.e., when the signal is extremely ``sparse'').
For instance, comparing the three examples mentioned above to
Bonferroni asymptotically as $m\rightarrow \infty$, we see that the corresponding step-up procedures
have a rejection function larger than Bonferroni's threshold
--- and are therefore {\em a posteriori}
more  powerful than Bonferroni --- provided their number of rejections $\abs{R}$ is larger than:
\begin{itemize}
\item $\Theta(\log m)$ for $\nu(k) \propto k^{-1}$ (Benjamini-Yekutieli procedure);
\item $\Theta(\sqrt{m})$ for $\nu$ uniform;
\item $\Theta\paren{\paren{m}^{\frac{2}{3}}}$ for $\nu(k) \propto k$\,.
\end{itemize}
(Recalling here that $\Theta()$ means asymptotic order of magnitude, in other terms ``asymptotically lower and upper bounded, up to a constant factor''.)
In each of the above cases, the largest proportion $u_m=\abs{R}/m$ of rejections for which Bonferroni's procedure
would {\em a posteriori} have been more powerful tends to zero as $m\rightarrow \infty$\,.
An identical conclusion will hold if we compare these rejection functions
to that of Holm's step-down \citep{Holm1979}, since the latter is equivalent to Bonferroni when $u_m \rightarrow 0$\,
(in addition, Holm's procedure is step-down while the above procedures are step-up).

More generally, let us exhibit a generic family of shape functions
$\beta$ such that $u_m$ tends to zero as $m\rightarrow \infty$. We
first define the proportion $u_m$ for a given shape function $\beta$
more formally, as $u_m=r_m/m$, where $r_m$ is the first point of
$\{1,\ldots,m\}$ for which $\beta(\cdot)$ is above $1$ (Bonferroni's
shape function). Introduce the family of \textit{scale
invariant} shape functions $\beta$, that is, the $\beta$s that can be
rewritten under the form $\beta(r)=m\wt{\beta}(\frac{r}{m})$ for some
\textit{fixed} function $\wt{\beta}(u)=\int_0^u v d\wt{\nu}(v)$ and
\textit{fixed} probability measure $\wt{\nu}$ on $(0,1]$. In the
latter, $\wt{\nu}$ should be taken independently of $m$ as a ``prior''
on the proportion of rejections. (Equivalently, $\wt{\nu}$ takes the role
of $\nu$ if we consider the following alternate scaling of the standard $\Lambda$-$\pi$
weighting: $\Lambda$ is the uniform probability measure on $\cH$ and $\pi \equiv 1$.)  It is then straightforward to check
that $u_m$ tends to $0$ as $m\rightarrow \infty$ if we choose
$\wt{\nu}$ such that $\wt{\beta}(u)>0$ for all $u>0$ (i.e. the origin
is an accumulation point of the support of $\wt{\nu}$).  This gives
many examples of shape functions which outperform Bonferroni's and
Holm's procedures as $m$ grows to infinity in the ``non-sparse''
case.
For example, the ``power function'' choice $d\wt{\nu}(u) = \ind{u \in [0,1]}
(\gamma+1) u^\gamma dx$ for $\gamma > -1$ gives rise to the rescaled
shape function $\wt{\beta}(u) = \frac{\gamma+1}{\gamma+2}
u^{\gamma+2}$ and thus
$\beta(r)=\frac{\gamma+1}{\gamma+2}\frac{r^{\gamma+2}}{m^{\gamma+1}}$. In
the cases $\gamma=0,1$, note that the latter corresponds to the
functions $\beta$ considered earlier (up to discretization).

By contrast, one can easily check that there is no scale-invariant {\em linear} rejection
function satisfying \eqref{condbeta_chapself}:
the Benjamini-Yekutieli procedure would correspond (up to lower order terms introduced by discretization)
to the ``truncated'' prior
$d\wt{\nu}(u)=(\log m)^{-1}\ind{m^{-1}\leq u \leq 1} x^{-1}du$\,, which cannot be extended to the origin
independently of $m$ since $u\mapsto u^{-1}$ is not integrable in 0.
We have seen above that $u_m \rightarrow 0$ nevertheless also holds
for this procedure: hence scale-invariant shape functions are certainly
not the only candidates in the family to
asymptotically outperform Bonferroni's and Holm's procedures in the ``non-sparse''
case.\looseness=-1

For comparison w.r.t. several other possible choices of $\nu$\,, (and for a finite $m=1000$) we have
systematically added Holm's rejection function on the plots of Figure \ref{fig_choicenu}. This
leads to a qualitatively similar conclusion.

\subsection{Adaptive step-up procedures}\label{appli_adapt}

We now give a very simple application of our results in the framework of {\em adaptive} step-up procedures.
Observe that the FDR control obtained for classical step-up
procedures is in fact not at the target level $\alpha$\,, but rather
at the level $\pi_0 \alpha$\,, where $\pi_0=\Pi(\cH_0)$ is the
``weighted volume'' of the set of true null hypotheses (equal to the
proportion of true null hypotheses $m_0/m$ in the standard case). This
motivates the idea of first estimating $\pi_0^{-1}$\, from the data
using some estimator $G(\bp)$\,, then applying the step-up procedure
with the modified shape function $\wt{\beta} = G(\bp)\beta$\,. Because
this function is now data-dependent, establishing FDR control for the
resulting procedure is more delicate; it is the subject of numerous
recent works (see, e.g., \citealp{Black2004,BKY2006,FDR2008}; see also \citealp{GBS2008} for an adaptive step-down procedure).

In this context we prove the following simple result, which is valid under the different types of
dependency conditions:
\begin{lemma}
\label{lem_dep}
Assume either of the following conditions is satisfied:
\begin{itemize}
\item the $p$-values $(p_h,h\in\cH)$ are PRDS on $\cH_0$\,, $\beta$ is the identity function.
\item the $p$-values have unspecified dependencies and $\beta$ is a function of the form~\eqref{condbeta_chapself}.
\end{itemize}
Define $R$ as an adaptive step-up procedure using the data-dependent threshold collection
 $\Delta(h,r,\bp) = \alpha_1 \pi(h) G(\bp)\beta(r)$\,, where $G(\bp)$\, is some estimator of  $\pi_0^{-1}$\,,
assumed to be nondecreasing as a function of the $p$-values.
Then the following inequality holds:
\begin{equation}
\FDR(R) \leq \alpha_1 + \e{\frac{|R\cap \cH_0|}{|R|} \ind{ G(\bp) > \pi_0^{-1} }}\,. \label{lemma_adap_df}
\end{equation}
\end{lemma}
\begin{proof}
Consider $\wt{R}$ the modified step-up procedure using the data-dependent threshold collection $\alpha_1\pi(h)\max(\pi_0^{-1},G(\bp))\beta(r)$\,.
Then it is easy to check that $\wt{R}$ satisfies the self-consistency condition {\bf SC}$(\alpha_1\pi_0^{-1},\pi,\beta)$.
Furthermore, $\wt{R}$ is a nondecreasing set as a function of the $p$-values, by the hypothesis
on $G$\,. Therefore, by combining Proposition~\ref{mlem} with Proposition~\ref{cor_PRDSimpliesfdrcontrol}
(resp. Proposition~\ref{lemmaCCdep}), $\wt{R}$ has its FDR controlled at level $\pi_0(\alpha_1\pi_0^{-1})=\alpha_1$
in both dependency situations and we have
\begin{align*}
\FDR(R) & = \e{\frac{|R\cap \cH_0|}{|R|}\ind{|R| > 0 }}\\
&  \leq \e{\frac{|\wt{R}\cap \cH_0|}{|\wt{R}|}\ind{|\wt{R}| > 0 }} + \e{\frac{|R\cap \cH_0|}{|R|} \ind{ G(\bp) > \pi_0^{-1} }}\\
& \leq \alpha_1 + \e{\frac{|R\cap \cH_0|}{|R|} \ind{ G(\bp) > \pi_0^{-1} }}\,.
\end{align*}
\end{proof}
Incidentally, the above proof illustrates a technical use of the main result where the
inclusion in the self-consistency condition is generally not an equality.

We can apply  Lemma~\ref{lem_dep} when considering a so-called {\em two-stage}
procedure, where $\pi_0$ is estimated using a preliminary multiple
testing procedure $R_0$\,. We assume here that this first stage has controlled
FWER (e.g. Holm's step-down).
\begin{corollary}
Let $R_0$ be a multiple testing procedure with
$FWER(R_0):=\P(\cH_0\cap R_0\neq
\emptyset)$ controlled at level $\alpha_0$\,.
Estimate $\pi_0$ by
$\wh{\pi}_0=\Pi((R_0)^{c})=\sum_{h\notin R_0}\pi(h)\Lambda(\set{h})$ the
$\pi$-volume of hypotheses non rejected by the first stage, and put
$G(\bp)=\wh{\pi}_0^{-1}$ (defined as $+\infty$ when $\wh{\pi}_0=0$).

Then the adaptive step-up procedure $R$ using the data-dependent threshold collection
 $\Delta(h,r,\bp) = \alpha_1 \pi(h) G(\bp)\beta(r)$ satisfies
\[
\FDR(R) \leq \alpha_0 + \alpha_1\,.
\]
\end{corollary}
The proof is a direct application of Lemma \ref{lem_dep}:
the second term in \eqref{lemma_adap_df} is upper bounded by $\P(G(\bp) > \pi_0^{-1})=\P(\Pi((R_0)^c)<\Pi(\cH_0))$, which is itself smaller than or equal to $\P(\cH_0\cap R_0\neq \emptyset)$, the FWER of the first stage.
Note that in the standard situation where $\Lambda=|.|$ is the
counting measure and $\pi$ is uniform, the above estimator of
$\pi_0^{-1}=m/m_0$ is simply $m/\wh{m}_0$, where $\wh{m}_0$ is the
number of non rejected hypotheses by the first stage.

Because of the loss in the level introduced by the first
stage, the latter result is admittedly not extremely sharp: for example, putting
$\alpha_0=\alpha_1 = \alpha/2$\,, a theoretical improvement over the
non-adaptive version at level $\alpha$ is obtained only when more than
50\% of hypotheses are rejected in the first stage. However, while
sharper results are available under the assumption of independent
$p$-values \citep[see, e.g.,][]{BKY2006}, up to our knowledge, there are
almost no results addressing the case of dependent $p$-values (as is
the case in the above result). The results we know of for this case
are found in works of \citet{Sar2006} and \citet{Far2007}. The latter reference
establishes a result similar to the above one, but seems to make the
implicit assumption that the two stages are independent, which we are
not assuming here. A more extensive treatment of the question of
adaptive procedures when following the general principles exposed in
the present work, including other applications of Lemma \ref{lem_dep},
is proposed by \citet{BR2008b} (see also the second author's PhD thesis, \citealp{Roq2007}, Chap. 11).

\subsection{FDR control over a continuous space of hypotheses}\label{sec:continuoushypo}

An interesting feature of the approach advocated here for proving FDR control is
that it can be readily adapted to the case where $\cH$ is a continuous
set of hypotheses. A simple example where this situation arises
theoretically is when
when the underlying observation is modelled as a random process $W$ over a
continuous space $\cT$\,, and the goal is to test for each $t \in \cT$
whether $\e{W(t)}=0$\,. In this case we can identify $\cH$ to $\cT$\,.
Such a setting was considered for example by
\citet{PGVW2004}.

In order to avoid straying too far from our main message in the present work, it was
decided to postpone the detailed exposition of this point to a
separate note. We refer the interested reader to the Section 5 of
the technical report of \citet{BR2008}, and restrict ourselves here to a brief overview.
First, under appropriate (and tame)
measurability assumptions, the framework developed in this paper
carries over without change: in the FDR definition, instead of using
the cardinality measure (which is of course not adapted in the
continuous case), we are able to deal with an arbitrary ``volume
measure'' $\Lambda$ on $\cH$ (such as the Lebesgue measure if $\cH$ is a compact subset of $\mbr^d$).
Also, while it seems considerably more
difficult to define rigorously step-up procedures in the traditional
sense via reordering of the $p$-values, Definition
\ref{def_stepup_chapself} of a step-up procedure carries over
in a continuous setting.

Secondly, our main tool, Proposition~\ref{mlem}, remains true when
$\cH$ is continuous, by replacing each sum over $\cH$ by the corresponding
integral (with respect to the measure $\Lambda$). Thirdly comes the question
of how to adapt the three types of dependency conditions considered in
Section \ref{sec:depcond} to a continuous setting. Under unspecified dependencies,
there is nothing to change as our arguments are not specific to the discrete setting.
The independent case, on the other hand, cannot be adapted to the continuous setting
as it conflicts with some measurability assumptions. However, this setting is
mainly irrelevant in a continuous setting as continuous families
of independent random variables are not usually considered.
Finally, in the case of positive dependencies,
condition \eqref{R-nondec} still ensures the dependency control
condition since Lemma \ref{UVlemma} is valid for arbitrary variables,
not necessarily discrete.  The main difficulty is therefore to
suitably adapt the PRDS assumption in the continuous setting.
We propose two extensions of the PRDS condition,
namely the ``strong continuous PRDS'', which is a direct
adaptation of the finite PRDS definition to a continuous setting,
and the ``weak continuous PRDS'', which states that any finite
subfamily of $p$-values should be (finite) PRDS. The strong continuous
PRDS condition is sufficient but arguably possibly not easy to check,
while the weak PRDS condition is easier but requires some additional
requirements on the procedure $R$ to ensure condition {\bf DC}.
An example of a process satisfying either type of condition is a continuous Gaussian
process with a positive covariance operator.

\vspace{-3pt}
\subsection{Other types of procedures}
\label{ssec:nonsu}
\vspace{-3pt}

We want to point out that the approach advocated here
also provides FDR control for procedures more general than step-up.
For example, as mentioned at the end of Section \ref{section_stepup},
generalized step-up-down procedures satisfy a self-consistency
property. Therefore, combining Proposition~\ref{mlem} with
Proposition~\ref{cor_PRDSimpliesfdrcontrol} (PRDS case) and
Proposition~\ref{lemmaCCdep} (unspecified dependencies), we obtain the
following result:

\begin{theorem}\label{theostep-wise}
Assume either of the following conditions is satisfied:
\begin{itemize}
\item the $p$-values $(p_h,h\in\cH)$ are PRDS on $\cH_0$\,, $\beta$ is the identity function.
\item the $p$-values have unspecified dependencies and $\beta$ is of the form \eqref{condbeta_chapself}.
\end{itemize}
Then the generalized step-up-down procedure of any order $\lambda\in [0,|\cH|]$ and associated
to the threshold collection $\Delta(h,r) = \alpha \pi(h) \beta(r)$ has its FDR controlled at level $\alpha \Pi(\cH_0)$\,.
\end{theorem}
In the PRDS case and with the standard $\Gamma$ - $\pi$ weighting, the first point of the above result has been first proved by \citet{Sar2002}
(see also \citealp{FDR2008}, where an approach related to ours is used to prove the same result; this is discussed in more detail
below in Section \ref{sec:discuss}).
The contribution of the above result is to deal with possible $\Gamma$ - $\pi$ weighting and
with the general dependent case (in particular, note that this theorem contains both Theorem~\ref{th_stepup} and Theorem~\ref{cor_marteau}).
We emphasize that the latter result does not come trivially from the
  fact that a step-up-down procedure is always a subset of the step-up procedure
using the same threshold collection, because in the FDR expression the
numerator and the denominator inside the expectation {\em both} decrease
with the rejection set size.

It could  however legitimately be objected that only step-up procedures are
really of interest in the present context, since they are less conservative
than step-up-down, and even the less conservative possible under the {\bf
SC} condition, as argued in Section \ref{section_stepup}. But one interest of the self-consistency
condition is to allow more flexibility, in particular if there are {\em additional}
constraints to be taken into account. Consider the following plausible scenario: in a medical imaging
context, the user wants to enforce additional geometrical constraints on the
set $R$ of rejected hypotheses, represented as a 2D set of pixels. For example,
one could demand that $R$ be convex or  have only a limited number of connected
components. If such additional constraints come into play, the step-up may not
be admissible, and has to be replaced by a subset satisfying the constraints. In this case,
the flexibility introduced by the {\bf SC} condition will be useful in order to
give a simple criterion sufficient to establish FDR control
without necessarily having to engineer a new proof for each new specific algorithm.
Note in particular that in such a scenario, one would probably like to choose a maximal rejection set
satisfying both the geometric constraints and self-consistency
condition; in this case the resulting procedure cannot be characterized in general as a step-up-down procedure,
and the {\bf SC} condition might hold without equality, i.e. $R\subsetneq L_{\Delta}(|R|)$\,.

\subsection[Another application of condition DC$(\beta)$]{Another application of condition \ref{DCC}}
\label{ssec:stepdown}

In this section, we step outside of the framework used in Proposition \ref{mlem};
more precisely, we present another application of condition \ref{DCC} to study
the FDR of a step-down procedure that does not satisfy the self-consistency condition with respect to the adequate shape function.
We will prove that the step-down procedure proposed by
\citet{BL2} and \citet{RS2006} has a controlled FDR under a PRDS-type
assumption of $\cH_0$ on $\cH_1$; we also deduce
a straightforward generalization to the unspecified dependencies case.
In this section, we only consider $\Lambda$ equal to the counting measure,
so that the aim is to control the standard FDR.

\citet{BL2} and \citet{RS2006} introduced the step-down procedure
based on the threshold collection $\Delta(i)=\frac{\alpha m}{(m-i+1)^2}$\,,
showed that it has controlled
FDR at level $\alpha$ if for each $h_0\in\cH_0$\,, $p_{h_0}$ is independent of the collection of
$p$-values $(p_{h},h\in\cH_1)$ (in fact \citealp{RS2006} used a
slightly weaker assumption, but it reduces to independence when the $p$-values
of true null hypotheses are uniform on $[0,1]$).  Here, we prove this result under
a weaker assumption, namely
a positive regression depency assumption of $p$-values of $\cH_1$ from those of $\cH_0$\,.
Let us reformulate slightly the notion of ``PRDS on
$\cH_0$'' given in Definition~\ref{def_PRDS}.
We say that the $p$-values of
$(p_h,h\in\cH_1)$ are positively regression dependent from each one
in a separate set $\cH_0$ (for short: $\cH_1$ PRDSS on $\cH_0$)
if for any measurable nondecreasing set $D\subset[0,1]^{\cH_1}$ and for all
$h_0\in\cH_0$, the function $$u\mapsto
\Proba\big(({p_{h}})_{h \in\cH_1}\in D \telque p_{h_0} =u\big)$$ is
nondecreasing. Note that the latter condition is obviously satisfied
when for all $h_0 \in \cH_0$\,, $p_{h_0}$ is independent of $(p_{h},h\in\cH_1)$\,.
We chose to introduce a new acronym only to emphasize the fact that, contrarily
to the standard PRDS\,, this assumption does not put constraints on the inner
dependency structure of the $p$-value vector of true hypotheses.

\begin{theorem}\label{th_stepdown}
Suppose that the $p$-values of $\cH_1$ are PRDSS on $\cH_0$.
Then the step-down procedure of threshold collection
$\Delta(i)=\frac{\alpha m}{(m-i+1)^2}$ has a FDR less than or equal to
$\alpha$\,.

If $\beta$ is a shape function of the form \eqref{condbeta_chapself}, then
without any assumptions on the dependency of the $p$-values, the step-down
procedure of threshold collection $\Delta(i) = \frac{\alpha m}{m-i+1} \beta\paren{\frac{1}{m-i+1}}$
has a FDR less than or equal to $\alpha$\,.
\end{theorem}
The proof is found in appendix. Essentially, we followed the proof of Benjamini and Liu
(\citeyear{BL2})
and identified the point where the condition \ref{DCC} (along with the results of Lemma \ref{UVlemma})
can be used instead of their argument.

\citet{BL1999} proposed a slightly less conservative step-down
procedure: the step-down procedure with the threshold collection
$\Delta(i)=1-\big[1-\min\big(1,\frac{\alpha
m}{m-i+1}\big)\big]^{1/(m-i+1)}.$ It was proved by \citet{BL1999} that this
procedure controls the FDR at level $\alpha$ as soon as the $p$-values
are independent. More recently, a proof of this result was given by
\citet{Sar2002} when the $p$-values are MTP$_2$ (see the definition
there) and if the $p$-values corresponding to true null hypotheses are
exchangeable. However, the latter conditions are more restrictive than
the PRDSS assumption of Theorem \ref{th_stepdown}.

The procedure of Theorem \ref{th_stepdown}
is often more conservative than the LSU procedure. First because the
LSU procedure is a step-up procedure, and secondly because the
threshold collection of the LSU procedure is larger on a substantial range.
However, in some specific cases ($m$ small and large number of
rejections), the threshold collection of Theorem \ref{th_stepdown} can
be larger than the one of the LSU procedure.
A similar argument can be made when comparing the proposed modified step-down
under unspecified dependencies to (for example) the modified LSU procedure
of \citet{BY2001}.

In order to use Theorem~\ref{th_stepdown} in the unspecified dependencies case, we have to choose a ``prior'' $\nu$ on the
set $\set{\frac{1}{k}: 1 \leq k \leq m}$\,:
\begin{itemize}
 \item taking a uniform $\nu$ yields 
 $\Delta(i) =  \alpha \frac{1}{m-i+1} \paren{ \frac{1}{m-i+1} + \cdots + \frac{1}{m}}\,,$
\item taking $\nu\paren{\frac{1}{k}} \propto k$ results in the threshold function $\Delta(i) = \frac{\alpha}{m+1} \frac{2i}{ (m-i+1)}\,,$
\item taking $\nu\paren{\frac{1}{k}} \propto \frac{1}{k}$ results in $\Delta(i)$ equal to \\
\noindent$\gamma_m^{-1} \alpha \frac{m}{m-i+1} \paren{ \frac{1}{(m-i+1)^2} + \cdots + \frac{1}{m^2}} \simeq \gamma_m^{-1} \alpha \frac{i}{(m-i+1)^2}\,,$
with $\gamma_m = \sum_{i\leq m} \frac{1}{i}$\,. 
\end{itemize}

\section{Discussion and conclusion}\label{concl_selfbound}

\subsection{The self-consistency condition and connection with other works}\label{sec:discuss}

The self-consistency condition with a linear shape function can be related to the following
heuristic motivation:
{consider the problem of choosing a threshold for rejected $p$-values, which we
reformulate equivalently as choosing $r$ such that $L_{\Delta}(r)$
has a FDR smaller than $\alpha$ (for the linear threshold collection $\Delta(h,r)=\alpha r/m$).
If the final number of rejections $|L_{\Delta}(r)|$ was equal to a {\em deterministic} constant $C(r)$, we
would have a FDR bounded by
\[
\e{\,\abs{\{h\in\cH_0\telque p_h \leq \alpha
r/m \}}\,}/ C(r) \leq \alpha r/ C(r)\,,
\]
so that the desired FDR control would be attained if $r\leq
C(r)=|L_\Delta(r)|$, that is, when $L_\Delta(r)$ satisfies the
self-consistency condition. This reasoning is, of course, unrigorous since
$L_{\Delta}(r)$ is in fact a random variable (and we need other arguments to correctly prove the FDR
control, e.g. Lemma~3.2). This point of view is
in the same spirit as the {\em
post-hoc} interpretation of the classical linear step-up procedure
proposed in Section~3.3 of
\citet{BH1995}, where the authors remarked that the
linear step-up procedure maximizes the number of rejected hypotheses under
the above constraint, which is the property we used in
Definition~\ref{def_stepup_chapself}.

As mentioned in the introduction and in Section \ref{ssec:nonsu}, the forthcoming paper of \citet{FDR2008}
introduces a condition quite similar to the self-consistency condition (although formulated differently).
Precisely, condition (T2) of\break \citet{FDR2008} can be seen to be equivalent to $R=L_{\Delta}(|R|)$
in our notation (in the specific case of a linear threshold collection $\Delta$ and for
the standard $\Lambda$-$\pi$ weighting).
It is proved in Theorem 4.1 of \citet{FDR2008} that (T2) implies
FDR control in the PRDS case (or more precisely, when \eqref{R-nondec} holds).
The authors note that the corresponding proof unifies and simplifies
classical results and proofs.  The present work, developed
independently, led to a very similar conclusion.  In particular, \citet{FDR2008} note that their result covers in general
the step-up-down procedures satisfying (T2), which is essentially the
same as the first point of the present Theorem~\ref{theostep-wise}
(for the standard $\Lambda$ and $\pi$-weighting).

As an additional contribution, we introduced the ``abstract'' dependency
condition {\bf DC}, which allowed us to
increase the range where the self-consistency condition
can be used, in particular when the $p$-values have unspecified
dependencies. We also included  $\Lambda$ and $\pi$-weighting
in our results; the formulation we adopted allows in particular for an
easy extension to infinite, possibly continuous hypothesis spaces.
Other original applications were exposed in Section~\ref{sec:appli}.

Conversely, \citet{FDR2008} used their approach for different
applications of interest,
based on an asymptotically optimal rejection curve. Several step-up or
step-up-down procedures are proposed by \citet{FDR2008} based on variations on this
rejection curve and shown to have a
an asymptotic and adaptive (in the sense of Section \ref{appli_adapt}) control of the FDR
(related to this is  also the step-down procedure of \citet{GBS2008}, based on the same
curve and shown to enjoy non-asymptotic control of the FDR). These
results do not fit directly into the framework delineated in the present paper,
but some of the technical tools used in their proof are of a similar spirit.
A full technical development on this topic is out of the scope of the present work,
but we demonstrate in a separate work \citep{BR2008b} that the two conditions
we presented here (along with some additional key ideas coming from \citealp{BKY2006})
can be used to prove (non-asymptotic) FDR control under independence,
for an adaptive procedure based on a rejection curve analogous to that considered
by \citet{FDR2008} and \citet{GBS2008}. To this regard, let
us also mention the recent work of \citet{Neu08}, which compares a number of these related
procedures in terms of their asymptotical power.

Finally, we mention that the self-consistency condition presented here has a slightly
weaker form than condition (T2) of  \citet{FDR2008}, namely it is $R\subset L_{\Delta}(|R|)$
instead of $R=L_{\Delta}(|R|)$.
From a technical point of view, we note here that
the argument of \citet{FDR2008} can actually be adapted straightforwardly to accomodate the weaker
condition. Is the weaker form of the condition of interest at all?
While the stricter condition is sufficient to cover the case of step-up and
step-up-down procedures, in the present work we have also tried to demonstrate that
the weaker form is not purely anectodical but useful in some other
applications:
first for truncated threshold collections (proof of Lemma~\ref{lem_dep}),
and secondly in Section \ref{ssec:nonsu} where we mentioned plausible practical scenarios where
equality might not hold due to additional constraints.

\subsection{Conclusion}

The approach advocated in this paper to establish FDR bounds
introduced a clear distinction between two sufficient
conditions of a different nature: on the one hand, the
self-consistency condition, which is purely algorithmic, and on the
other hand, the (essentially probabilistic) dependency control condition.
The two conditions are effectively coupled via the common
choice of the shape function $\beta$\, appearing in both. The
fundamental result of this paper is that these two conditions suffice
for FDR control, but part of our message is that this
point of view also introduced some relevant technical tools, which,
abstracting some key arguments present in previous works, can
be of use in various other settings.

While these conditions are only sufficient and hence certainly not universal,
we illustrated their interest by  recovering in Sections \ref{applistep-up}, \ref{ssec:unspec} and \ref{ssec:nonsu}
several existing results of the FDR multiple testing
literature in an unified way, as in particular with any arbitrary combination
of the following factors:
\begin{itemize}
\item arbitrary $\Lambda$-weighting of the FDR via the volume measure,
\item arbitrary $\pi$-weighting of the $p$-values via the weight function,
\item arbitrary choice of dependency setting: independent, PRDS or unspecified,
\item in the unspecified dependencies setting, arbitrary choice of the shape function $\beta$ satisfying \eqref{condbeta_chapself}.
\item in the procedure algorithm, arbitrary choice between ``step-down'' and ``step-up'', ``step-up-down'', and more
generally arbitrary choice among the possible orders $\lambda$ in a ``step-up-down'' procedure.
\end{itemize}
In the past literature, many results have been established for specific combinations of the above
variations; here we were able to cover all of these at once, possibly covering
combinations that had not been explicitly considered earlier (in particular, the
fourth ``factor'' above seems to be new). Several other applications were proposed.

An interesting direction for future work is
to try to ``adapt'' the choice of the weight function $\pi$ (and possibly also
the distribution $\nu$ in the case of unknown dependencies) depending on the observed data.
Because these parameters have an crucial influence on power, doing so in a principled way
might result in a substantial improvement.

\bigskip

\appendix

\noindent{\bf \Large Appendix}

\section{\texorpdfstring{Proof of Lemma \ref{UVlemma}}{Appendix A: Proof of Lemma 3.2}}

{\em Part (i).}
We want to establish the following inequality:
\[
\e{\frac{\ind{U\leq cg(U)}}{g(U)}} \leq c\,,
\]
for $U$ stochastically lower bounded by a uniform distribution and $g$
nonincreasing.
Let $\mathcal{U} = \{u \telque cg(u) \geq u\}$\,, $u^* = \sup\mathcal{U}$ and $C^{*}=\inf\{g(u)\telque u \in \mathcal{U}\}$\,.
It is not difficult to check that $u^* \leq c C^*$\, (for instance take any nondecreasing sequence
$u_n \in \mathcal{U} \nearrow u^*$\,, so that $g(u_n) \searrow C^*$\,).
If $C^*=0$\,, then $u^*=0$ and the result is trivial.
Otherwise, we have
\[
\e{\frac{\ind{U\leq cg(U)}}{g(U)}} \leq \frac{\Proba(U \in
 \mathcal{U})}{C^*} \leq \frac{\Proba( U \leq u^*) }{C^*} \leq \frac{u^*}{C^*} \leq
c.
\]

{\em Part (ii).}
The proof uses a similar telescopic sum argument as developed by \citet{BY2001}
for proving FDR control of the linear step-up under the PRDS assumption; the goal of the lemma
presented here is to isolate this argument in order to specifically concentrate on
condition {\bf DC}, and to extend it to arbitrary (non-discrete) variables\,.

We want to prove the inequality
\begin{equation}
\e{\frac{\ind{U\leq cV}}{V}} \leq c\,\nonumber
\end{equation}
for $U,V$ two nonnegative real variables such that
$U$ is stochastically lower bounded by a uniform
  distribution, and the conditional distribution of $V$ given $U\leq
  u$  is stochastically decreasing in $u$\,.
Fix some $\eps>0$ and some $\rho \in (0,1)$\, and choose $K$ large enough so that $\rho^K<\eps$.
Put $v_0 = 0$ and $v_i = \rho^{K+1-i}$ for $1\leq i \leq 2K+1$\,.
The following chain of inequalities holds:
\begin{multline*}
\e{\frac{\ind{U\leq cV}}{V \vee \eps }}\\
\begin{aligned} &
\leq \ \sum_{i=1}^{2K+1} \frac{\Proba(U\leq c v_i; V \in [v_{i-1}, v_i))}{v_{i-1} \vee \eps } + \eps\\
& \leq  c \sum_{i=1}^{2K+1} \frac{\Proba(U\leq c v_i ;
V \in [v_{i-1},v_i) )}{\Proba(U\leq c v_i)} \frac{v_i}{v_{i-1} \vee \eps} + \eps\\
& \leq c \rho^{-1} \sum_{i=1}^{2K+1} \Proba(V \in [v_{i-1},v_i) \telque U \leq c v_i) + \eps\\
& = c \rho^{-1} \sum_{i=1}^{2K+1} \big( \Proba(V < v_i\telque U \leq c v_i) - \Proba(V < v_{i-1}\telque U \leq c v_i) \big) + \eps\\
& \leq c \rho^{-1} \sum_{i=1}^{2K+1} \big( \Proba(V< v_i\telque U \leq c v_i) - \Proba(V< v_{i-1}\telque U \leq c v_{i-1}) \big) + \eps\\
& \leq c \rho^{-1} + \eps \,.
\end{aligned}
\end{multline*}
We obtain the conclusion by letting $\rho \rightarrow 1$\,, $\eps \rightarrow 0$ and applying the monotone
convergence theorem.

{\em Part (iii).}
Rewriting for any $z>0$\,, $1/z=\int_{0}^{+\infty}v^{-2}\ind{v\geq z}dv$\,, and using Fubini's theorem:
  \begin{align*}
\e{\frac{\ind{U\leq c\beta(V)}}{V}}&= \e{\int_{0}^{+\infty} v^{-2} \ind{v\geq V} \ind{U\leq c\beta(V)} dv }\\
&=\int_{0}^{+\infty} v^{-2}\e{\ind{v\geq V} \ind{U\leq c \beta(V)}}dv\\
&\leq \int_{0}^{+\infty} v^{-2}\Proba\big(U\leq c\beta(v)\big)dv\\
&\leq c\int_{0}^{+\infty} v^{-2}\beta(v)dv\\
&= c \int_{u\geq 0} u \int_{v \geq 0 } \ind{u \leq v} v^{-2} dv d\nu(u) =c\,.
\end{align*}
\qed

\section{\texorpdfstring{Proof of Theorem \ref{th_stepdown}}{Appendix B: Proof of Theorem 4.6}}\label{sec_appendix}

To establish the first assertion of the Theorem, remember we assume the $p$-values of $\cH_1$ are PRDSS on $\cH_0$\,,
and the threshold collection is $\Delta(i)=\alpha m/(m-i+1)^2$\,.
Assume $m_0>0$ (otherwise the result is trivial) and consider $p_{(1)}\leq p_{(2)} \leq \dots \leq p_{(m)}$ the ordered
$p$-values of $(p_h,h\in\cH)$. Denote by $j_0$ the
(data-dependent) smallest integer $j\geq 1$ for which $p_{(j)}$
corresponds to a true null hypothesis. Denote by $R_1$ the step-down
procedure of threshold collection $\Delta$ and restricted to the set
of the false null hypotheses $\cH_1$. First note that the
following points hold:
\begin{itemize}
\item[$(i)$] $|R\cap \cH_0|>0\Rightarrow  p_{(j_0)}\leq \frac{\alpha m}{(m-j_0+1)^2}$
\item[$(ii)$] $|R\cap \cH_0|>0 \Rightarrow j_0-1\leq |R_1|$
\item[$(iii)$] $R_1 \subset R \cap \cH_1$\,.
\end{itemize}
 To prove this, suppose that $|R\cap \cH_0|>0$, so that the null
hypothesis corresponding to $p_{(j_0)}$ is rejected by $R$. Hence,
from the definition of a step-down procedure we have $p_{(j_0)}\leq
\Delta(j_0)$ and $(i)$ holds. Moreover, since for all $j\leq j_0-1,$
we have $p_{(j)}\leq \Delta(j)$ and $p_{(j)}$ corresponds to a false
null hypothesis, $R_1$ necessarily rejects all the null hypotheses
corresponding to $p_{(j)},j\leq j_0-1$, and we get $(ii)$\,. Finally,
we obviously have $R_1 \subset \cH_1$ and it is easy to check that $R_1 \subset R$
(using the fact that the reordered $p$-values of $\cH_1$ form a subsequence of
$(p_{(i)})$).

From $(i)$ and $(ii)$ we deduce that
\begin{equation}
\label{eq_ext_ap}
\abs{R \cap \cH_0} > 0 \Rightarrow \exists h \in \cH_0\,: p_h \leq \frac{\alpha m}{(m-|R_1|)^2} \leq \frac{\alpha m}{m_0 (m-|R_1|)}\,.
\end{equation}
Therefore,
\begin{eqnarray}
\FDR(R)&=&\e{\frac{|R\cap \cH_0|}{|R|} \ind{|R\cap \cH_0|>0}}\nonumber\\[3pt]
& = & \e{\frac{|R\cap \cH_0|}{|R\cap \cH_0| + |R\cap \cH_1|} \ind{|R\cap \cH_0|>0}}\nonumber\\[3pt]
&\leq& \e{\frac{m_0}{m_0+|R\cap \cH_1|} \ind{|R\cap \cH_0|>0}}\nonumber\\[3pt]
&\leq& \sum_{h\in\cH_0}\e{\frac{m_0}{m_0+|R_1|}\ind{p_h \leq (\alpha m/m_0) (m-|R_1|)^{-1}}},\nonumber
\end{eqnarray}
where for the first inequality, we used that fact that for each fixed
$a\geq 0$, $x\mapsto \frac{x}{x+a}$ is a nondecreasing function on
$\R^+\backslash \{0\}$\,.  For the second inequality, we used
simultaneously \eqref{eq_ext_ap} and  the point $(iii)$ above.  Since the function
$x\mapsto\frac{m_0}{m_0+x}\frac{m}{m-x}$ is log-convex on $[0,m_1]$
and takes values $1$ in $x=0$ and $x=m_1$, we have pointwise
$\frac{m_0}{m_0+|R_1|}\frac{m}{m-|R_1|}\leq 1\,.$ Therefore, we
get\vadjust{\eject}
\begin{eqnarray}
\FDR(R)&\leq& \frac{1}{m}\sum_{h\in\cH_0}\e{\frac{\ind{p_h \leq (\alpha m/m_0) (m-|R_1|)^{-1}}}{(m-|R_1|)^{-1}}}\nonumber\\[3pt]
&\leq& \frac{1}{m}\sum_{h\in\cH_0} \alpha m/m_0 = \alpha \,\nonumber.
\end{eqnarray}
In the last inequality, we used that the couple $(p_h,(m-|R_1|)^{-1})$ satisfies condition \ref{DCC}
with $c=\alpha
m/m_0$ and $\beta(x)=x$\,; this holds
in the present case from part (ii) of Lemma~\ref{UVlemma} because for any $v>0$,
$D=\{\mathbf{z}\in [0,1]^{\cH_1}\telque
(m-|R_1(\mathbf{z})|)^{-1}<v\}$ is a nondecreasing set (so that we can
apply the same reasoning as for the proof of Proposition~\ref{cor_PRDSimpliesfdrcontrol}).

For the second part of the theorem, we follow exactly the same proof as above with the
modified threshold function and part (iii) of Lemma~\ref{UVlemma} instead of part~(ii).
\qed


\section*{Acknowledgements}

The authors wish to thank the anonymous reviewer and AE for their constructive criticism of previous versions of the manuscript,
which allowed to improve the overall organization and focus of the paper.

\bibliographystyle{apalike} 

\end{document}